\documentclass[12pt]{article}

\usepackage{amsmath}
\usepackage{amsfonts}
\usepackage{amssymb}
\usepackage{graphicx, color}
\usepackage{hyperref}
\usepackage{epsfig}
\usepackage{epstopdf} 

\definecolor{darkgreen}{rgb}{0,0.55,0}

\newtheorem{proposition}{Proposition}[section]
\newtheorem{theorem}{Theorem}[section]
\newtheorem{lemma}[theorem]{Lemma}

\newtheorem{remark}[theorem]{Remark}

\newtheorem{definition}{Definition}

\def\phi{{\varphi}}

\DeclareSymbolFont{AMSb}{U}{msb}{m}{n}
\DeclareMathSymbol{\N}{\mathbin}{AMSb}{"4E}
\DeclareMathSymbol{\Z}{\mathbin}{AMSb}{"5A}
\DeclareMathSymbol{\R}{\mathbin}{AMSb}{"52}
\DeclareMathSymbol{\Q}{\mathbin}{AMSb}{"51}
\DeclareMathSymbol{\I}{\mathbin}{AMSb}{"49}
\DeclareMathSymbol{\C}{\mathbin}{AMSb}{"43}

\begin{document}
\title{Least Gradient Problems with Neumann Boundary Condition}

\author{ Amir Moradifam\footnote{Department of Mathematics, University of California, Riverside, CA, USA. E-mail: moradifam@math.ucr.edu.} \qquad
}

\date{\today}

\smallbreak \maketitle

\begin{abstract}
We study existence of minimizers of the least gradient problem 
\[\inf_{v \in BV_g} \int_{\Omega}\varphi(x, Dv),\]
where $BV_g=\{v \in BV(\Omega): \int_{\partial \Omega}gv=1\}$, $\varphi(x,p):  \Omega\times \R^n \rightarrow \R$ is a convex, continuous, and homogeneous function of degree $1$ with
respect to the $p$ variable, and  $g$ satisfies the comparability condition $\int_{\partial \Omega} g dS=0$. We prove that  for every $0\not \equiv g \in L^{\infty}(\partial \Omega)$ there are infinitely many minimizers in $BV(\Omega)$. Moreover there exists a divergence free vector field $T\in (L^{\infty}(\Omega))^n$ that determines the structure of level sets of all minimizers, i.e. $T$ determines $\frac{Du}{|Du|}$, $|Du|-$ a.e. in $\Omega$, for every minimizer $u$. We also prove some existence results for general 1-Laplacian type equations with Neumann boundary condition. A numerical algorithm is presented that simultaneously finds $T$ and a minimizer of the above least gradient problem. Applications of the results in conductivity imaging are discussed. 

\end{abstract}

\section{Introduction and Statement of the Main Results}
Let $\Omega$ be a bounded open set in $\R^n$ with Lipschitz boundary and $\varphi: \Omega\times \R^n \rightarrow \R$ be a continuous function satisfying the following conditions:\\

($C_1$) There exists $\alpha_1, \alpha_2>0$ such that $\alpha_1|p| \leq \varphi(x,p) \leq \alpha_2 |p|$ for all $x\in \Omega$ and $p \in \R^n$. \\

($C_2$) $p \mapsto \varphi(x,p)$ is a norm for every $x$.\\ \\
This work is a continuation of the author's work on existence, uniqueness, and structure of minimizers of the least gradient problems in \cite{JMN, M-Indiana, MNTCalc}. In this paper we study the general least gradient problem 
\begin{equation}\label{LTVProb}
\inf_{v\in \mathcal{M}_g} \int_{\Omega} \varphi(x, Dv),
\end{equation}
where  $g \in L^{\infty}(\Omega)$ satisfies the compatibility condition 
\begin{equation}\label{compat}
\int_{\partial \Omega} g dS=0,
\end{equation}
and 
\begin{equation*}
\mathcal{\mathcal{M}}_g:=\{v  \in BV(\Omega)\ \  \hbox{and} \ \ \int_{\partial \Omega} g v dS=1 \}.
\end{equation*}
Such problems arise in conductivity imaging (see \S 1.1) and are closely related to the 1-Laplacian type equation 
\begin{eqnarray}\label{MainPDE}
\left\{ \begin{array}{ll}
\nabla_x \cdot \nabla_p \varphi(x,\frac{Du}{|Du|})=0  &\text{in } \Omega\\ \\
\left[ \nabla_p \varphi (x,\frac{Du}{|Du|}), \nu_{\Omega} \right]=\lambda g&\text{on }\partial \Omega,
\end{array} \right.
\end{eqnarray}
where $\lambda>0$ is a constant and $\frac{Du}{|Du|}$ is the Radon-Nikodym derivative of $Du$ with respect to $|Du|$, and the boundary condition is understood in the sense of the integration by parts formula \eqref{IBP-Trace} below. When $\varphi(x,p)=a(x)|p|$  for some positive function $a\in C(\bar{\Omega})$, then \eqref{LTVProb} reduces to the weighted least gradient problem
\[ \inf_{v \in \mathcal{M}_g} \int_{\Omega}a|Dv|,\]
and \eqref{MainPDE} reduces to the 1-laplacian equation 

\begin{eqnarray}\label{MainPDE1}
\left\{ \begin{array}{ll}
\nabla \cdot ( a \frac{Du}{|Du|})=0  &\text{in } \Omega\\ \\
\left[  a \frac{Du}{|Du|}, \nu_{\Omega} \right]=\lambda g&\text{on }\partial \Omega.
\end{array} \right.
\end{eqnarray}
Least gradient and 1-laplacian problems with Dirichlet boundary condition have been studied extensively in \cite{JMN, MRD, Mercaldo1, Mercaldo1, M-Indiana, sternberg_ziemer92, sternberg_ziemer, sternbergZiemer93}, and with Neumann boundary condition in \cite{Andr1, Mercaldo}. 

For any $v \in BV_{loc}(\R^n)$, let $\varphi(x,Dv)$ be the measure defined by
\begin{equation}\label{varphi.def}
\int_U \varphi(x,Dv) \ = \ \int_U \varphi(x, \frac{Dv}{|Dv|}) |Dv|
\qquad\mbox{ for any bounded Borel set $U$}.
\end{equation}
Standard  facts about $BV$ functions imply that (see \cite{AB}) if $\Omega$ is an open set,
then
\begin{equation}\label{varphi.def1}
\int_{\Omega} \varphi(x,Dv) = \sup \{ \int_{\Omega} v \nabla \cdot X dx\ \ : \ \ X \in C^{\infty}_{c} (\Omega; \R^n), \ \ \sup \varphi^0(x, X(x)) \leq 1 \}, 
\end{equation}
where $\varphi^0(x, \cdot)$ is defined by
\begin{equation}\label{varphi0}
\varphi^0(x,\xi) =\sup \{\frac{\xi \cdot p}{\varphi(x,p)}: \ \ p\in \R^n\}, 
\end{equation}
(see \cite{AB, JMN, M-Indiana}). 
For $v\in BV(\Omega)$, $\int_{\Omega}\varphi(x, Dv)$ is called the $\varphi$-total variation of $v$ in $\Omega$. 
Let $\nu_{\Omega}$ denote the outer unit normal vector to $\partial \Omega$. For every $V\in (L^{\infty}(\Omega))^n$ with div$(V) \in L^n(\Omega)$ there exists a unique function $[V,  \nu_{\Omega}] \in L^{\infty}(\partial \Omega)$ such that
\begin{equation}\label{IBP-Trace}
\int_{\partial \Omega} [V,  \nu_{\Omega}] u d\mathcal{H}^{n-1}=\int_{\Omega} u \nabla \cdot V dx +\int_{\Omega}  V \cdot \nabla u dx, \ \ \forall u \in C^1(\bar{\Omega}).
\end{equation}
In addition, for $u\in BV(\Omega)$ and $V\in (L^{\infty}(\Omega))^n$ with div$(V) \in L^n(\Omega)$, the linear functional $u\mapsto(V \cdot Du)$ gives rise to a Radon measure on $\Omega$, and
\begin{equation}\label{IBP0}
\int_{\partial \Omega} [V, \nu_{\Omega}]u d\mathcal{H}^{n-1}=\int_{\Omega} u \nabla \cdot V dx +\int_{\Omega}  (V \cdot D u), \ \ \forall u \in BV(\Omega),
\end{equation}
see \cite{Al, An, And4}. We are now ready to give a precise definition of solutions of \eqref{MainPDE} (see Definition 4.2 in \cite{Mercaldo}).

\begin{definition}\label{solution}
We say that $u\in BV(\Omega)$ is a solution to \eqref{MainPDE} if there exists a vector field $T$ such that 
\begin{equation}\label{test1}
\nabla \cdot T = 0 \hbox{ on } D'(\Omega),
\end{equation}
\begin{equation}\label{test2}
\nabla_p \varphi(x,\frac{Du}{|Du|})=T, \ \ |Du|-\hbox{a.e.} \ \ \hbox{in} \ \ \Omega,
\end{equation}
\begin{equation}\label{test3}
\left[T,\nu_{\Omega}\right]  =  \lambda g, \ \ \mathcal{H}^{n-1}-\hbox{a.e. on} \ \ \partial \Omega.
\end{equation}
\end{definition}

\begin{definition}\label{solution0}
A function $u \in BV(\Omega)$ is said to be an entropy solution of \eqref{MainPDE} if it is a solution to \eqref{MainPDE} in the sense of Definition \ref{solution}, and
\begin{equation}\label{test4}
\varphi(x,p) \geq  T \cdot p, \ \ \forall p\in \R^n \ \ \hbox{and a.e. $x$ in} \ \ \Omega. 
\end{equation}
\end{definition}
We shall prove that minimizers of the least gradient problem \eqref{LTVProb} are entropy solutions of the 1-Laplacian type equation \eqref{MainPDE}. 

\begin{remark} Consider the special case $\varphi(x,p)=a(x)|p|$ where $a(x) \in C(\overline{\Omega})$ is a positive function. If $u\in BV(\Omega)$ is a solution of the 1-laplacian equation \eqref{MainPDE1} in the sense of Definition \ref{solution0}, then 
\[a \frac{Du}{|Du|}=T\]
for some divergence free vector field $T$ with $|T|\leq a$ a.e. and $[T,\nu_{\Omega}]=g$ ($\lambda=1$). For a physical interpretation of this definition, assume $J$ is the current density induced inside a conductive body $\Omega$ with electrical conductivity $\sigma \in C(\Omega)$. Then by Ohm's law $\nabla \cdot J \equiv 0$ and
\[J=-\sigma Du.\]
Hence if we let  $a=|J|$, then the corresponding voltage potential $u$ satisfies the equation \eqref{MainPDE1} and the induced current density vector field $J$ plays the role of $T$ in Definitions \ref{solution} and \ref{solution0} (see \S 1.1 below). 

\end{remark}

The following theorem settles the question of existence of minimizers of \eqref{LTVProb}. 

\begin{theorem}\label{theoLG}
Suppose $\varphi: \Omega \times \R^n \rightarrow \R$ is a continuous function satisfying the condition $(C_1)$ and $(C_2)$, and $0 \not \equiv g\in L^{\infty}(\partial \Omega)$ satisfies the compatibility condition \eqref{compat}. Then the least gradient problem \eqref{LTVProb} admits infinitely many minimizers in $\mathcal{M}_g$. Moreover,  there exists a vector field $T \in (L^{\infty}(\Omega))^n$ satisfying \eqref{test1}, \eqref{test3}, and \eqref{test4} with
\[\lambda:= \inf _{v \in \mathcal{M}_g}\int_{\Omega} \varphi(x,Dv)\]
such that 
\begin{equation}\label{directionInside}
\varphi(x, \frac{Du}{|Du|})= T\cdot \frac{Du}{|Du|}, \ \ |Du|-\hbox{a.e.} \ \ \hbox{in} \ \ \Omega,
\end{equation}
for every minimizer $u$ of \eqref{LTVProb}. In particular all minimizers of \eqref{LTVProb} have the same level set structure and are entropy solutions of the 1-Laplacian equation \eqref{MainPDE}. 
\end{theorem}

\begin{remark}
The above theorem asserts that a fixed divergence free vector field $T$ determines the structure of the level sets of all minimizers of  the least gradient problem \eqref{LTVProb}. More precisely, since 
\[\varphi(x,p)\geq T\cdot p  \]
for every $p\in S^{n-1}$ and a.e.  $x\in  \Omega$, it follows from \eqref{directionInside} that $|Du|$-a.e., $p=\frac{Du}{|Du|}$ maximizes 
\[\frac{T \cdot p}{\varphi(x,p)}\]
among all $p\in S^{n-1}$, determining  $\frac{Du}{|Du|}$, $|Du|$-a.e. in $\Omega$. This is a remarkable fact about minimizers of least gradient problem \eqref{LTVProb}. In the special case $\varphi(x,p)=a|p|$, Theorem \ref{theoLG} implies that for every minimizer $u$ of \eqref{LTVProb} 
\[\frac{Du}{|Du|}\cdot T=|T|=a \ \ |Du|-a.e. \ \ \hbox{in} \ \ \Omega,\]
i.e.  $\frac{Du}{|Du|}$ is parallel to $T$, $|Du|-a.e.$ . Similar phenomenon occurs for minimizers of general least gradient ptoblems with Dirichlet boundary condition \cite{M-Indiana}. 
\end{remark}

The following proposition describes the connection between the solutions of \eqref{MainPDE} and minimizers of the least gradient problem \eqref{LTVProb}. 

\begin{proposition}\label{propRelation}
Suppose $\varphi: \Omega \times \R^n \rightarrow \R$ is a continuous function satisfying the condition $(C_1)$ and $(C_2)$, and $0 \not \equiv g\in L^{\infty}(\partial \Omega)$ satisfies the compatibility condition \eqref{compat}. Then $u\in\mathcal{M}_g$ is a minimizer of \eqref{LTVProb} if and only if it is an entropy solution of \eqref{MainPDE} in the  sense of Definitions \ref{solution} and \ref{solution0} with 
\begin{equation}
\lambda:= \inf _{v \in \mathcal{M}_g}\int_{\Omega} \varphi(x,Dv).
\end{equation}
\end{proposition}

\begin{proposition} \label{propMultiple}
Let $u$ be a solution of \eqref{MainPDE} and $F$ be an increasing Lipschitz continuous function. Then $F(u)$ is also a solution of \eqref{MainPDE}. \\
\end{proposition}
The next results follows immediately from Theorem \ref{theoLG}, Proposition \ref{propRelation}, and Proposition \ref{propMultiple}. 

\begin{theorem}\label{theoPDE}
Suppose $\varphi: \Omega \times \R^n \rightarrow \R$ is a continuous function satisfying the condition $(C_1)$ and $(C_2)$ and $0 \not \equiv g\in L^{\infty}(\partial \Omega)$ satisfies the compatibility condition \eqref{compat}. Then there exists $\lambda^*>0$ such that the equation (\ref{MainPDE}) with $\lambda=\lambda^*$ has infinitely many entropy solutions in $\mathcal{M}_g$, and for $\lambda\neq \lambda^*$  (\ref{MainPDE}) does not admit any entropy solutions. Moreover, when $\lambda=\lambda^*$, there exists $T \in (L^{\infty}(\Omega))^n$ satisfying \eqref{test1}, \eqref{test2}, and \eqref{test3} such that 
\[\varphi(x, \frac{Du}{|Du|})= T\cdot \frac{Du}{|Du|}, \ \ |Du|-\hbox{a.e.} \ \ \hbox{in} \ \ \Omega,\]
for every solution $u$ of \eqref{MainPDE}. In particular all solutions of \eqref{MainPDE} have the same level set structure. 
\end{theorem}

In \cite{Mercaldo} authors studied the existence of solutions of \eqref{MainPDE} for the special case $\varphi(x,p)=|p|$ by analyzing the behavior of solutions of the p-laplacian problem 
\begin{eqnarray}\label{p-laplace}
\left\{ \begin{array}{ll}
\nabla \cdot \left( |\nabla u|^{p-2}\nabla u \right)=0  &\text{in } \Omega\\ \\
 |\nabla u|^{p-2} \frac{\partial u}{\partial \nu_{\Omega}}= g&\text{on }\partial \Omega,
\end{array} \right.
\end{eqnarray}
as $p\rightarrow 1$, and showed that if $||g||_*=1$, then solutions of \eqref{p-laplace} converge to a solution of
\begin{eqnarray}\label{1-laplace}
\left\{ \begin{array}{ll}
\nabla \cdot \left( \frac{\nabla u}{|\nabla u|} \right)=0  &\text{in } \Omega\\ \\
\left[ \frac{Du}{|Du|}, \nu_{\Omega} \right]= g&\text{on }\partial \Omega,
\end{array} \right.
\end{eqnarray}
where 
\[||g||_{*}:=\sup_{v\in \mathcal{S}_1\setminus \{0\}} \left\lbrace \frac{\int_{\partial \Omega} g u dS}{\int_{\Omega}|\nabla u|dx}   \right\rbrace, \]
and
\[\mathcal{S}_1=\left\lbrace v\in W^{1,1}(\Omega): \int_{\partial \Omega} v dS=0 \right\rbrace.\]
Note that $||g||_*=1$ corresponds to the case $\lambda =1$ in \eqref{MainPDE}. If $||g||_* < 1$ or $||g||_* >1 $, then solutions of \eqref{p-laplace} converge to $u\equiv 0$, or $\infty$ on a set of positive measure, respectively  \cite{Mercaldo}. Therefore if $||g||_* \neq 1$, then solutions of \eqref{MainPDE} can not be obtained as a limit of solutions of \eqref{p-laplace} without the knowledge of the parameter $\lambda$. Moreover, the convergence is extremely unstable with respect to perturbations of $\lambda$. In Section 3, we shall present a numerical algorithm for solving \eqref{MainPDE} which simultaneously finds $\lambda^*$, $T$, and a solution of the degenerate equation (\ref{MainPDE}). This algorithm converges to a solution of (\ref{MainPDE}) with $\lambda=\lambda^*$ independent of the value of $||g||_*$.

\subsection{Applications in Conductivity Imaging}
Least gradient problem \eqref{LTVProb} arise in  the inverse problem of determining an electrical conductivity $\sigma$ of a conductive body $\Omega$ from one measurement of the magnitude of the current density field $|J|$ generated inside $\Omega$ and the voltage potential $f$ on the boundary $\partial \Omega$.  Indeed if the electrical conductivity is isotropic, then the voltage potential inside $\Omega$ is the unique minimizer of the least gradient problem 
\[\inf_{u\in BV_f(\Omega)} \int_{\Omega}a|D u|dx,\]
where $a=|J|$ and $BV_f(\Omega):=\{u\in BV(\Omega): u|_{\partial \Omega}=f\}$ (see \cite{MNT, MNTa_SIAM, NTT07, NTT08, NTT10, NTT11}). Once $u$ is determined inside $\Omega$, then the conductivity $\sigma$ can be easily determined inside $\Omega$.

One can also consider the inverse problem of recovering an isotropic conductivity $\sigma$ from the knowledge of the magnitude of the induced current $a=|J|$ and $g=J\cdot \nu_{\Omega}$ on $\partial \Omega$. Notice that 
\[\nabla \cdot J=\nabla \cdot (a\frac{Du}{|Du|})=0,\]
and hence the voltage potential is a minimizer of the least gradient problem 
\[\inf_{u\in \mathcal{M}_g} \int_{\Omega}a|D u|dx,\]
or equivalently $u$ is a solution of \eqref{MainPDE1}. It follows from Theorem \ref{theoLG} that the voltage potential $u$ and consequently the conductivity $\sigma$ can not be uniquely identified from the knowledge of $|J|$ inside $\Omega$ and $g=J\cdot \nu_{\Omega}$ on $\partial \Omega$. However, the full current density vector field $J$ can be uniquely recovered $a(x) d\mu$-a.e. ($d\mu$ is the Lebesgue measure). The current density vector field $J$ is indeed the vector field $T$ in Definitions \ref{solution} and \ref{solution0} which is also a solution of the corresponding dual problem described in Section 2 below.

In \cite{HMN} the author and his collaborators presented a method for recovering the conformal factor of an anisotropic conductivity matrix in a known conformal class from one interior measurement of current density.  Suppose the matrix valued conductivity $\sigma(x)$ is of the form
\[\sigma(x) =c(x)\sigma_0(x)\]
where $c(x)\in C^{\alpha}(\Omega)$ is a positive scalar valued function and $\sigma_0 \in C^{\alpha}(\Omega, Mat(n, \R^n))$ is  a known positive definite symmetric matrix valued function. The conformal factor  $\sigma_0$ can be determined using Diffusion Tensor Magnetic Resonance Imaging. They showed that the corresponding voltage potential $u$ is the unique minimizer of the least gradient problem 
\[\inf_{u \in BV_f(\Omega)} \int_{\Omega} \varphi(x, Dv),\]
where $\varphi$ is given by 
\begin{equation}\label{varphi}
\varphi(x,p)=a(x)\left( \sum_{i,j=1}^{n}\sigma_0^{ij}(x)p _i p_j\right)^{1/2},
\end{equation}
\begin{equation}\label{aJJ}
a=\sqrt{\sigma_0^{-1} J \cdot J},
\end{equation}
and $J$ is the induced current density vector field. One may also similarly consider the problem of recovering a current density vector field $J$ induced by an anisotropic conductivity $\sigma$ from the knowledge of the conformal factor $\sigma_0$, $a=\sqrt{\sigma_0^{-1} J \cdot J}$, and $g=J \cdot \nu_{\Omega}$ on $\partial \Omega$. Then the corresponding voltage potential will be a minimizer the least gradient problem \eqref{LTVProb} where $\varphi$ is given by \eqref{varphi}. Similar to the isotropic case, $u$ can not be uniquely recovered. However, by Theorem \ref{theoLG}, the current density vector field $J$ can be uniquely determined $|Du|$-a.e. in $\Omega$, where $u$ is an arbitrary solution of \eqref{LTVProb}.  See also \cite{NTV} where weighted least gradient problems are utilized to analyze conductivity imaging form the knowledge of the magnitude of the induced current density vector field with complete electrode model boundary conditions.

\section{Proof of the Main Results} 
Let $\Omega$ be a bounded open set in $\R^n$ with Lipschitz boundary, and $0 \not \equiv g\in L^{\infty}(\partial \Omega)$. Choose $u_g \in W^{1,1}(\Omega)$ with $\int_{\partial \Omega}g u_g dS=1$. Define 
\[\mathcal{W}^{0}_g(\Omega):= \{u \in W^{1,1}(\Omega): \ \ \int_{\partial \Omega}ug=0 \},\]
and
\[\mathcal{W}^{1}_g(\Omega):= \{u \in W^{1,1}(\Omega): \ \ \int_{\partial \Omega}ug=1 \}.\]
Let $E: (L^1(\Omega))^n\rightarrow \R$ and $G: \mathcal{W}^{0}_g(\Omega) \rightarrow \R$ be defined as follows
\begin{equation}\label{FandG}
E(P):=\int_{\Omega} \varphi (x, P+\nabla u_g)dx, \ \ G(u)\equiv 0.
\end{equation}
Then the problem (\ref{LTVProb}) can be written as 
\[(P) \ \ \ \ \ \inf_{u \in \mathcal{W}^{0}_g(\Omega)} E(\nabla u)+G(u).\]
By Fenchel duality (see Chapter III in \cite{ET}) the dual problem is given by  
\begin{equation*}
(P^*) \hspace{.5cm}\sup_{V \in (L^{\infty}(\Omega))^n} \{ -E^* (-V)-G^*(\nabla^* V)\}.
\end{equation*}
where $E^*$ and $G^*$ are the Legendre-Fenchel transform of $E$ and $G$ respectively, and $\nabla^*$ is the adjoint of $\nabla :\mathcal{W}^{0}_g(\Omega) \rightarrow (L^{1}(\Omega))^n$. Recall that the Legendre-Fenchel transform  $E^*: (L^{\infty}(\Omega))^n \rightarrow \R$ is defined as follows 
\[E^*(V)=\sup  \{\langle V, P \rangle -E(P): \ \  P \in  (L^1(\Omega))^n\}.\]
The following lemma is proved by the author in \cite{M-Indiana}. 

\begin{lemma}[\cite{M-Indiana}] \label{f^*lem} Let  $E$ be defined as in equation (\ref{FandG}). Then 

\begin{equation}\label{f^*}
E^*(V)=\left\{ \begin{array}{ll}
-\langle D u_g, V\rangle\ \ &\hbox{if} \ \ \varphi^0(x,V(x))\leq 1 \ \  \hbox{in} \ \ \Omega,\\
\infty  \ \ & \hbox{otherwise}.
\end{array} \right.
\end{equation}
\end{lemma}

\vspace{.2cm}

\begin{lemma} Let $G: \mathcal{W}^{0}_g(\Omega) \rightarrow \R$ be defined as $G \equiv 0$. Then for $G^*: (\mathcal{W}^{0}_g (\Omega))^*\rightarrow \R$ we have

\begin{equation}G^*(\nabla^* b)=\left\{ \begin{array}{ll}
0 &\hbox{if} \ \ b \in \mathcal{B}\\
\infty   &\hbox{otherwise},
\end{array} \right.
\end{equation}\label{G*} 
where 
\[\mathcal{B}:=\{b \in (L^{\infty}(\Omega))^n:  \nabla \cdot b \equiv 0 \hbox{ and } b\cdot \nu_{\Omega}=\lambda g\ \ \mathcal{H}^{n-1}-a.e. \hbox{ on } \partial \Omega, \hbox{ for some } \lambda \in \R\}.\]

\end{lemma}
{\bf Proof.} By definition 
\begin{eqnarray*}
G^*(\nabla^* b)=\sup_{u\in \mathcal{W}^0_g} \langle \nabla^* b, u\rangle= \sup_{u\in \mathcal{W}^0_g} \langle b, \nabla u\rangle=\sup_{u\in \mathcal{W}^0_g} \left(  \int_{\partial \Omega} [b,\nu_{\Omega}]u dS-\int_{\Omega} u \nabla \cdot b dx \right).
\end{eqnarray*}
Since $W^1_0 \subset \mathcal{W}^0_g$, the above supremum will be $\infty$ if $\nabla \cdot b \not \equiv 0$. Hence we have 
\begin{equation*}\label{FunctionalAnalysis}
G^*(\nabla^* b)=\sup_{u\in \mathcal{W}^0_g}  \int_{\partial \Omega} [b,\nu_{\Omega}]u dS
\end{equation*}
Also since $\mathcal{W}^0_g$ is a vector space, the above supremum will be infinity unless 
\begin{equation}\label{FunctionalAnalysis}
\int_{\partial \Omega} [b,\nu_{\Omega}]u dS=0 \ \ \hbox{ for all } u\in \mathcal{W}^0_g.
\end{equation}
Now let $N=\{\lambda g: \lambda \in \R\} \subset L^{\infty}(\partial \Omega).$ Then $\mathcal{W}^0_g|_{\partial \Omega}=N^{\perp}$, and it follows from the second geometric form of the Hahn-Banach theorem (see Proposition 1.9 in \cite{Brezis} for a proof) that 
\[(N^{\perp})^{\perp}=\bar{N}=N.\]
Hence $[b,\nu_{\Omega}]\in N$ and the proof is complete. \hfill $\Box$ \\

Define
\[\mathcal{V}:=\{V\in \mathcal{B}:\ \ \varphi^0(x,V(x))\leq 1 \ \ \hbox{in} \ \ \Omega \}.\]
It follows from Lemmas \ref{f^*lem} and \ref{G*} that the dual problem can be explicitly written as
\begin{equation*}\label{expDual} 
(P^*) \ \ \ \ \ \ \ \ \sup_{V \in \mathcal{V}}  \int _{\partial \Omega} [V, \nu_{\Omega}] u_g dS,  
\end{equation*}
where $\nu_{\Omega}$ is outward pointing unit normal vector on $\partial \Omega$.  The primal problem (P) may not have a solution, but the dual problem ($P^*$) always has a solution. This is a direct consequence of Theorem III.4.1 in \cite{ET}. Indeed it easily follows from \eqref{varphi.def1} that $I(v)=\int_{\Omega}\varphi (x, Dv)$ is convex, and $J: L^1(\Omega)\rightarrow \R$ with $J(p)=\int_{\Omega} \varphi(x, p)dx$ is continuous at $p=0$ (a consequence of $C_2$). Therefore the condition (4.8) in the statement of Theorem III.4.1 in \cite{ET} is satisfied and the following proposition holds.

\begin{proposition}\label{dualityGap}
Let $\Omega \subset \R^n$ be a bounded open set with Lipschitz boundary and assume $\varphi: \Omega \times \R^n \rightarrow \R$ be a continuous function satisfying the condition $(C_1)$ and $(C_2)$, and $0 \not \equiv g \in L^{\infty}(\partial \Omega)$ satisfies the compatibility condition \eqref{compat}.  Then there exists a divergence free vector field $T \in (L^{\infty}(\Omega))^n$ with $\varphi^0(x,T)\leq 1$ a.e. in $\Omega$ such that
\begin{eqnarray*}
\inf_{v\in \mathcal{M}_g} \int_{\Omega} \varphi(x, Dv) = \max_{V \in \mathcal{V}} \int _{\partial \Omega}  [V, \nu_{\Omega}] u_g dS= \int _{\partial \Omega}  [T, \nu_{\Omega}] u_g dS. 
\end{eqnarray*}
In particular the dual problem $P^*$ has a solution $T \in \mathcal{V}$.
\end{proposition}
Proof of Proposition \ref{propMultiple} is similar to the proof of Theorem 4.3 in \cite{Mercaldo} and we omit it. \\ \\
{\bf Proof of Proposition \ref{propRelation}.} Let $u$ be a minimizer of the least gradient problem (\ref{LTVProb}). By Proposition \ref{dualityGap}, the dual problem has a solution $T\in \mathcal{V}$ with $[T, \nu_{\Omega}]=\lambda g$ for some $\lambda \in \R$. Since $\varphi^0(x,T(x))\leq 1$ a.e. in $\Omega$,  \eqref{test4} holds and 
\begin{eqnarray*}
\int_{\Omega} \varphi(x, Du)&=& \int_{\Omega}\varphi(x,\frac{Du}{|Du|})|Du|\geq \int_{\Omega}T\cdot \frac{ Du}{|Du|}|Du|\\
&=& \int_{\Omega} T\cdot Du= \int_{\partial \Omega}u [T,\nu_{\Omega}] dS\\
&=&\int_{\partial \Omega}\lambda u g dS=\lambda \\
&=&  \int_{\partial \Omega} [T,\nu_{\Omega}]u_g dS=\int_{\Omega} \varphi(x, Du).
\end{eqnarray*}
Thus the inequality is indeed an equality. Hence 
\begin{equation}\label{equiv-gradient}
\varphi(x, \frac{Du}{|Du|})= T\cdot \frac{Du}{|Du|}, \ \ |Du|-\hbox{a.e.} \ \ \hbox{in} \ \ \Omega.
\end{equation}
Now consider the function 
\[K(x,p)=\varphi(x,p)-T(x)\cdot p.\]
It follows from \eqref{test4} and \eqref{equiv-gradient} that for $|Du|$-a.e. $x \in \Omega$, $K_x(p)=K(x,p)$ attains its minimum at $p=\frac{Du}{|Du|}$. Thus 
\[\nabla_p K (x,\frac{Du}{|Du|})=\nabla_p(x, \frac{Du}{|Du|})-T(x)=0, \ \ |Du|-a.e. \ \ x\in \Omega.\]
Hence \eqref{test2} holds and $u$ is an entropy solution of \eqref{MainPDE}. 

Conversely, assume $u\in \mathcal{M}_g$ is a entropy solution of \eqref{MainPDE}. Since $\varphi(x,p)$ is a homogeneous functions of order $1$ with respect to the $p$ variable,
\[\varphi(x,\frac{Du}{|Du|})=\nabla_p \varphi(x,\frac{Du}{|Du|}) \cdot \frac{Du}{|Du|}=T\cdot \frac{Du}{|Du|}, \ \ |Du|-a.e. \ \ x\in \Omega.\]
Therefore it follows from the above computations that $T$ is a solution of the dual problem $(P^*)$ and $u$ is a minimizer of the least gradient problem \eqref{LTVProb}. \hfill $\Box$ \\ \\ 
{\bf Proof of Theorem \ref{theoLG}.} Suppose 
\[\beta:=\inf_{v\in \mathcal{M}_g} \int_{\Omega} \varphi(x, Dv)=\inf_{v\in W^0_{g}}\int_{\Omega}\varphi(x,Du_g+Dv).\]
Let $v\in \mathcal{M}_g$. It follows from ($C_1$) and continuity of the trace operator that 
\begin{eqnarray*}
\beta=\int_{\partial \Omega} g v dS& \leq & ||g||_{L^{\infty}(\partial \Omega)}\int_{\partial \Omega}|v|dS \\
&\leq& C ||g||_{L^{\infty}(\partial \Omega)} \int_{\Omega}|Dv|\\
&\leq & \frac{C ||g||_{L^{\infty}(\partial \Omega)}}{\alpha_1} \int_{\Omega}\varphi(x,Dv).
\end{eqnarray*}
Thus $\beta>0$. Now let $\{u_n\}_{n=1}^{\infty} $ be a minimizing sequence in $\mathcal{M}_g$, i.e. 
\[\lim_{n\rightarrow \infty} \int_{ \Omega} \varphi(x,Du_n)dx=\beta. \]
Then there exists a subsequence of $\{u_{n_k}\}_{k=1}^{\infty}$ that converges weakly$^*$ in $BV(\Omega)$ to some $u\in BV(\Omega)$, i.e. $u_{n_k}\rightarrow u$ strongly in $L^2(\Omega)$ and $Du_{n_k} \rightharpoonup Du$ weakly in the sense of measures. Since $I(u)=\int_{\Omega}\varphi(x,Du)$ is weakly lower semicontinuous (see \cite{JMN}), 
\[\int_{\Omega} \varphi(x, Du) \leq \beta.\]
Now let $T$ be a solution of the dual problem $(P^*)$ whose existence is guaranteed by Proposition \ref{dualityGap}. Then $[T, \nu_{\Omega}]=\lambda g$ for some $\lambda \in \R$ and we have 

\begin{eqnarray*}
\int_{\Omega} \varphi(x, Du)&=& \int_{\Omega}\varphi(x,\frac{Du}{|Du|})|Du|\geq \int_{\Omega}T\cdot \frac{ Du}{|Du|}|Du|\\
&=& \int_{\Omega} T\cdot Du= \lim_{k \rightarrow \infty} \int_{\Omega} T\cdot Du_{n_k} \\
&=&\lim_{k \rightarrow \infty} \int_{\partial \Omega}u_{n_k} [T,\nu_{\Omega}] dS=\lim_{k \rightarrow \infty} \int_{\partial \Omega}\lambda u_{n_k} g dS=\lambda \\
&=&\int_{\partial \Omega}u_g [T,\nu_{\Omega}] dS=\beta.\\
\end{eqnarray*}
Therefore $\beta=\lambda > 0$,
\[\int_{\Omega} \varphi(x, Du)=\beta,\]
and
\begin{equation}\label{equT}
\varphi(x,\frac{Du}{|Du|})=T\cdot \frac{Du}{|Du|}, \ \ |Du|-\hbox{a.e. in} \ \ \Omega.
\end{equation}
Moreover 
\begin{eqnarray}
\int_{\partial \Omega} ug dS&=&\frac{1}{\lambda}\int_{\partial \Omega} u[T,\nu_{\Omega}]dS\\
&=& \frac{1}{\lambda} \int_{\Omega} T\cdot Du =1. 
\end{eqnarray}
Hence $u\in M_g$ is a minimizer of \eqref{MainPDE}. 

Let $F$ be an increasing Lipschitz continuous function with $\int_{\partial \Omega}F(u)g dS \neq 0$. Then there exists $c_1, c_2 \in \R$ such that $c_1F(u)+c_2 \in \mathcal{M}_g$. Thus by Proposition \ref{propMultiple} equation \eqref{MainPDE} admits infinitely entropy solutions satisfying \eqref{test1}-\eqref{test4} for a fixed vector field $T \in \mathcal{V}$. By Proposition \ref{propRelation} the least gradient problem \eqref{LTVProb} also has infinitely many minimizers in $\mathcal{M}_g$. \hfill $\Box$

\section{An Algorithm for Finding Solutions}
In this section we present a numerical algorithm for solving the equation (\ref{MainPDE}) or equivalently finding a minimizer of \eqref{LTVProb}. Since the equation \eqref{MainPDE} is degenerate and the least gradient problem \eqref{LTVProb} does not have a unique minimizer, developing a numerical algorithm for finding such minimizers is in general challenging. Assuming that \eqref{LTVProb} has a minimizer in $u \in H^1(\Omega)$, we develop an algorithm that generates two sequences $(u_k)_{k\geq 1}$ and $(b_k)_{k\geq 1}$ such that $u\rightharpoonup u$ weakly in $H^1(\Omega)$ and $b_k\rightharpoonup T$ weakly in $L^2(\Omega)$, where $u$ and $T$ are solutions of \eqref{LTVProb} and its dual problem ($P^*$), respectively. In applications to conductivity imaging, it is natural to assume that the conductivity $\sigma$ belongs to $L^{\infty}(\Omega)$, and hence the corresponding voltage potential $u$ belongs to $ H^1(\Omega)$. Therefore the algorithm we develop here can be applied to problems arising from conductivity imaging. We conjecture that even if \eqref{LTVProb} is only assumed to have a minimizer in $\mathcal{M}_g$, then the sequences$(u_k)_{k\geq 1}$ and $(b_k)_{k\geq 1}$ produced by our algorithm would still converge to a minimizer $u$ of \eqref{LTVProb}  weakly$^*$ in $BV(\Omega)$ and to a solution $T$ of ($P^*$) weakly in $L^{\infty}(\Omega)$, respectively.  

Suppose \eqref{LTVProb} has a minimizer in $H^1(\Omega)$ and $u_g \in H^1(\Omega)$ satisfies $\int_{\partial \Omega} u_g gdS=1$. Then (\ref{LTVProb}) can be written as 
\begin{equation}\label{oldP}
\hspace{1cm}\inf_{u \in \mathcal{H}_g} F(\nabla u)+G(u)
\end{equation}
where  
\[\mathcal{H}_g:=\{v \in H^1(\Omega): \ \ \int_{\partial \Omega} v dS=0\ \ \hbox{and } \int_{\partial \Omega} g v=0\},\]
and $F:(L^2(\Omega))^n\rightarrow \R$ and $G: \mathcal{H}_g \rightarrow \R$ are defined as follows
\begin{equation}\label{FandG}
E(d):=\int_{\Omega} \varphi(x,d+\nabla u_g), \ \ \hbox{and} \ \ G\equiv 0. 
\end{equation}
As described in Section 2, the dual problem can be written as 
\begin{equation} \label{dualPrime}
\sup_{V\in (L^{2}(\Omega))^n} \{-E^*(-V)-G^*(\nabla^* V) \}.
\end{equation}
Let us aim to find a minimizer $T$ of the dual problem (\ref{dualPrime}) which will determine the structure of the level sets of all minimizers of \eqref{FandG}. If $T$ is a minimizer of (\ref{dualPrime}),  then  
 \begin{equation}\label{inclusion}
 0\in A(T)+B(T),
 \end{equation}
 where $A:=\partial (G^*o(-\nabla ^*))$ and $B:=\partial F^*$ are maximal monotone set-valued operators on $(L^2(\Omega))^n$, since they are sub-gradient of convex, proper, lower semi-continuous functions (see \cite{FPA, Rb}). We will apply Douglas-Rachford splitting algorithm, described below, to solve (\ref{inclusion}). 

For a set-valued function $P: H\rightarrow 2^{H}$, let $J_P$ denote its resolvent i.e.,
\[J_{P}=(Id+P)^{-1}.\] 
Let $H$ be a real Hilbert space and $A, B: H\rightarrow 2^{H}$ be two set-valued maximal monotone operators. 
Note that if $P$ is maximal monotone, then the resolvent $J_P$ is single valued \cite{FPA, Rb}.  Lions and Mercier \cite{LM}  showed that for any general maximal monotone operators $A,B$ and any initial element $S_0$, the sequence defined by the Douglas-Rachford recursion:
\begin{equation}\label{iter:DR}
S_{k+1}=(J_A(2J_B-Id)+Id-J_B)S_k,
\end{equation}
converges weakly to some point $S \in \emph{H}$ such that $T=J_B (S)$ solves the inclusion problem (\ref{inclusion}).  Recent results also prove weak convergence of the sequence
$T_k=J_{ B}(S_k)$ to $T$ ( see \cite{S}, and Chapters 25 and 27 in \cite{BC}). The following theorem describes the Douglas-Rachford splitting  algorithm and summarizes the convergence results in \cite{LM, S}. 

\begin{theorem}\label{DRS}
Let $H$ be a Hilbert space and let $A, B: \emph{H}\rightarrow 2^{\emph{H}}$ be maximal monotone operators and assume that a solution of (\ref{inclusion}) exists. Then, for any initial elements $S_0$ and $T_0$ and any $\alpha>0$, the sequences $S_{k}$ and $T_k$ defined by
\begin{eqnarray}\label{DRS-27}
S_{k+1}&=&J_{\lambda A}(2T_k-S_k)+S_k-T_k \nonumber \\
T_{k+1}&=&J_{\alpha B}(S_{k+1}),
\end{eqnarray}
converge weakly to some $S$ and $T$ respectively. Furthermore, $T=J_{\alpha B}(S)$  and $T$ solves 
\[0\in A(T)+B(T).\]
\end{theorem}
To apply the Douglas-Rachford splitting algorithm to the operators $A:=\partial (G^*(\nabla ^*))$ and $B:=\partial F^*$, we need to evaluate the resolvents $J_{\alpha A}(2T_k-S_k)$ and $J_{\alpha B}(T_{k+1})$ at each iteration. The following lemma provides a method for computing such resolvents (see  \cite{Setzer, Setzer1} for a proof). 

\begin{lemma} \label{SetLem}
Let $H_1$ and $H_2$ be two Hilbert spaces, $f: H_1\rightarrow \R \cup \{\infty\}$ and a bounded linear operator $L:H_1\rightarrow H_2$. Assume that $\hat{v}$ is a solution of
\[\hat{v}=argmin_{v\in H_1} \{\frac{\alpha}{2}\parallel Lv+q\parallel^2+f(v)\}.\]
Then
\begin{equation}
 \alpha (L\hat{v}+q)=J_{\alpha \partial(f^*o(-L^*))}(\alpha q). \
\end{equation} 
\end{lemma}
\vspace{.9cm}

Given $S_k$ and $T_k$, let $u^{k}$ and $d^{k}$ be the minimizers of the functionals 
\[I_1(u)=\parallel \frac{2T_k-S_k}{\alpha} +\nabla u\parallel^2,\]
and
\[I_2(d)=F(d)+\frac{\lambda}{2}\parallel \frac{S_k}{\alpha} -d \parallel^2,\]
respectively. Then by Lemma \ref{SetLem} we have
\[J_{\alpha A}(2T_k-S_k)=\alpha \nabla u^{k}+2T_k-S_k,\]
and 
\[T_{k}=J_{\alpha B}(S_{k})=S_{k}-\alpha d^{k}.\]
From (\ref{DRS-27}) we have

\[S_{k+1}:=S_k-T_k+[2T_k-S_k+\alpha \nabla u^{k+1}]=S_k+\alpha [\nabla u^{k+1}-d^{k}].\]
Thus for $k \geq 1$ we have
\[S_{k}=S_0+\alpha \sum_{i=1}^{k-1} ( \nabla  u^{i}-d^{i}) +\alpha \nabla u^{k}, \ \ T_{k}=S_0+\alpha \sum_{i=1}^{k} (\nabla  u^{i}-d^{i}),\]
and
\[2T_k-S_k=S_0+\alpha \sum_{i=1}^{k} (\nabla  u^{i}-d^i) -\alpha d^{k}.\]
So if we let $b^k =\frac{S_0}{\alpha}+\sum_{i=0}^{k} (\nabla  u^{i}-d^{i})$, then 

\begin{equation}
S_k=\alpha (b^k+d^k), \ \ \ \ T_k=\alpha  b^k, \ \ k\geq 1. \\
\end{equation}
Therefore to evaluate $J_{\alpha A}(2T_k-S_k)$ and $J_{\alpha B}(S_{k+1})$ in (\ref{DRS-27}) for all $k\geq 0$,  it suffices to find the minimizers $u^{k+1}$ and $d^{k+1}$ of the functionals 
\begin{equation}\label{Iu1}
I^{k+1}_1(u)=\parallel \nabla u+b^k-d^k\parallel ^2
\end{equation}
on $\mathcal{H}_g$, and
\begin{equation}\label{Iu2}
I^{k+1}_2(d)=\int_{\Omega} a|d+\nabla u_g|dx+\frac{\alpha}{2}\parallel b^k+\nabla u^{k+1}-d\parallel^2,
\end{equation}
on $(L^2(\Omega))^n$, and set $b^{k+1}=b^k+\nabla u^{k+1}-d^{k+1}$. 

Minimizers of (\ref{Iu1}) in $\mathcal{H}_g$ satisfy the Euler-Lagrange equation 
\begin{equation}\label{gap}
\Delta u^{k+1}=\nabla \cdot (d^k-b^k)\ \ \hbox{with} \ \ \frac{\partial{u^{k+1}}}{\partial \nu}=\beta g+(d^k-b^k)\cdot \nu \ \ \hbox{on } \partial \Omega,
\end{equation}
for some $\beta \in \R$. Conversely, if $u \in \mathcal{H}_g$ is a solution of \eqref{gap} for some $\beta \in \R$, then $u$ is a minimizer of \eqref{Iu1}. To identify the parameter $\beta$ and find a minimizer of  \eqref{Iu1} in $\mathcal{H}_g$, let $w$ be a solution of  
\begin{equation}\label{gapw}
\Delta w=0\ \ \hbox{with} \ \ \frac{\partial{w}}{\partial \nu}= g \ \ \hbox{on } \partial \Omega.
\end{equation}
Since $g \not \equiv 0$, we have 
\begin{eqnarray}
0<\int_{\Omega}|\nabla w|^2dx=-\int_{\Omega}w \frac{\partial w}{\partial \nu}dS.
\end{eqnarray}
In particular, 
\[\int_{\Omega}wgdS\neq 0.\]
Now let $u^{k+1}$ be a solution of 
\begin{equation}\label{gap}
\Delta u^{k+1}=\nabla \cdot (d^k-b^k)\ \ \hbox{with} \ \ \frac{\partial{u^{k+1}}}{\partial \nu}=(d^k-b^k)\cdot \nu\ \ \hbox{on } \partial \Omega,
\end{equation}
and define
\[\beta^{k+1}=-\frac{\int_{\Omega}u^{k+1}gdS}{\int_{\Omega}wg}.\]
Then $v^{k+1}=u^{k+1}+\beta^{k+1} w$ is a minimizer of \eqref{Iu1} in $\mathcal{H}_g$. Note that this minimizer is unique up to adding a constant. 

On the other hand, in general, the minimizer of the functional $I^{k+1}_2(d)$ can be usually computed explicitly. For instance if $\varphi(x,p)=a|p|$, then

\[d^{k+1}(x)=\left\{ \begin{array}{ll}
\max \{|w^{k+1}(x)|-\frac{a}{\lambda},0\}\frac{w^{k+1}(x)}{|w^{k+1}(x)|} -\nabla u_g(x) &\hbox{ if }  |w^{k+1}(x)|\neq 0,\\
-\nabla u_g(x)   &\hbox{ if }  |w^{k+1}(x)|=0,
\end{array} \right.
\] 
where $w^{k+1}=\nabla u^{k+1}+\nabla u_g +b^k$. Hence we arrive at the following algorithm that simultaneously solves the problem (\ref{oldP}) and its dual problem \eqref{dualPrime}. \\ \\
 {\bf The Algorithm:}\\
 
Let $\alpha>0$, $u_g \in H^1(\Omega)$ with $\int_{\partial \Omega} g u_g dS=1 $, and initialize $b^0, d^0 \in (L^2(\Omega))^n$. Let $w$ be a solution of \eqref{gapw} with $\int_{\Omega}wdx=0$. For $k\geq 0$:\\ 
\begin{enumerate}
\item (a) Solve \[ \Delta u^{k+1}= \nabla \cdot (d^k(x)-b^k(x)), \ \  \frac{\partial u^{k+1}}{\partial \nu}=(d^k-b^k)\cdot \nu,\]
with $\int_{\Omega}u^{k+1}dx=0.$\\ 
(b) Compute 
\[\beta^{k+1}=-\frac{\int_{\Omega}u^{k+1}gdS}{\int_{\Omega}wg}\]
and set $v^{k+1}=u^{k+1}+\beta^{k+1} w$.
\item Compute $d^{k+1}$ by minimizing \eqref{Iu2}. 
\item Let
\[b^{k+1}(x)=b^{k}(x)+\nabla v^{k+1}(x)-d^{k+1}(x).\]
\end{enumerate}

The following theorem follows directly from  Theorem \ref{DRS} and guarantees convergence of the above algorithm.

\begin{theorem}\label{theoConvergegence}
Let $a \in L^2(\Omega)$ be a non-negative function and  $g\in L^{\infty}(\partial \Omega)$. Suppose \eqref{LTVProb} has a minimizer in $H^1(\Omega) \cap \mathcal{M}_g$. Then for any $b^0, d^0 \in (L^2(\Omega))^n$ the sequences $\{b^k\}_{k\in N}$, $\{d^k\}_{k\in N}$, and $\{v^k\}_{k\in N}$ produced by Algorithm 1 converge weakly (in $(L^2(\Omega))^n$, $(L^2(\Omega))^n$, and $H^1(\Omega)$, respectively) to some $\frac{T}{\alpha}$, $d$, and $v^*$. Moreover $v:=v^*+u_g$ is a solution of the minimization problem (\ref{LTVProb}),  $T$ is a solution of the dual problem \eqref{dualPrime}, and $d=\nabla v^*$. In other words $u$ and $T$ satisfy the conditions of Definition 1 with
\[\lambda:= \inf _{v \in \mathcal{M}_g}\int_{\Omega} \varphi(x,Dv).\]

\end{theorem}

Algorithm 1 is in the spirit of the alternatiing split Bregman algorithm proposed by Goldstein and Osher \cite{GO} in finite dimensional settings in image processing. As pointed out by Esser \cite{Es} and Setzer \cite{Setzer1}, the idea to minimize $I_1^{k+1}$ and $I_2^{k+1}$ alternatingly was first presented for the augmented Lagrangian algorithm by Gabay and Mercier \cite{GM} and Glowinski and Marroco \cite{GMar}. The resulting algorithm is called the alternating direction method of multipliers (ADMM) \cite{Gab} and is indeed equivalent to the alternating split Bregman algorithm. The convergence of ADMM and the alternatiing split Bregman algorithm in finite dimensional Hilbert spaces was established  by Eckstein and Bertsekas \cite{EcBe} and independently by Cai, Osher, and Shen \cite{COS} and  Setzer \cite{Setzer, Setzer1}. Motivated by least gradient problems arising in conductivity imaging in infinite dimensional Hilber spaces, the second author and his collaborator first proved convergence of the split Btegman type algorithms with Dirichlet boundary conditions in \cite{MNT}. \\

{\bf Acknowledgments.} The author is partially supported by an start-up grant from University of California at Riverside.

\end{document}